\documentclass[a5paper]{book}
\usepackage[russian]{babel}
\usepackage{amsthm,amsmath,amssymb,color,multirow,graphicx,subfigure,float}
\usepackage{inputenc}
\usepackage{wrapfig, float,array}
\newtheorem{thm}{Theorem}[section]

\theoremstyle{definition}
\newtheorem{defn}[thm]{Definition}
\newtheorem{rem}[thm]{Remark}

\numberwithin{equation}{section}

\def\d{\,{\rm d}}
\def\E{\,{\rm e}}

\newcommand{\be}{\begin{equation}}
\newcommand{\ee}{\end{equation}}

\begin{document}

\begin{flushright}
\begin{tabular}{l}
{\sf Uzbek  Mathematical}\\
{\sf Journal, 2019,\ No 1, pp.{29}-{42}}\\
\end{tabular}
\end{flushright}

\sloppy

\begin{center}
\textbf{\large On an optimal interpolation formula in $K_2(P_2)$ space}
\textbf{Babaev S.S.}
\end{center}

{\small\textbf{Abstract.} The paper is devoted to construction of an optimal interpolation formula in $K_2(P_2)$ Hilbert space.
Here the interpolation formula consists of a linear combination $\sum_{\beta=0}^NC_{\beta}(z)\varphi(x_\beta)$ of given values of a function $\varphi$ from the space $K_2(P_2)$. The difference between functions and the interpolation formula is considered as a linear functional called the error functional. The error of the interpolation formula is estimated by the norm of the error functional.
We obtain the optimal interpolation formula by minimizing the norm of the error functional by coefficients $C_{\beta}(z)$ of the interpolation formula. The obtained optimal interpolation formula is exact for trigonometric functions $\sin\omega x$ and $\cos\omega x$.
At the end of the paper we give some numerical results which confirm our theoretical results.

\textbf{Keywords:} extremal function, error functional, Hilbert space, optimal interpolation formula, optimal coefficients.

\textbf{Mathematics Subject Classification (2010):} 41A15.}

\makeatletter
\renewcommand{\@evenhead}{\vbox{\thepage \hfil {\it Babaev S.S.}\\[0.5mm]   \hrule }}
\renewcommand{\@oddhead}{\vbox{\hfill
{\it Optimal interpolation formulas in $K_2(P_2)$ space}\hfill
\thepage \hrule}} \makeatother

\label{Babaev1}

\setcounter{section}{0}

\section{Introduction and statement of the Problem}

There are algebraic and variational approaches in the spline theory. In the
algebraic approach splines are considered as some smooth piecewise polynomial
functions. In the variational approach splines are elements of Hilbert or Banach
spaces minimizing certain functionals. The first spline functions were constructed
from pieces of cubic polynomials. After that, this construction was modified, the
degree of polynomials increased. The theory of splines based on variational methods studied and
developed, for example, by J.H.Ahlberg et al. \cite{Ahlb67},    C. de Boor \cite{deBoor63}, A.R.Hayotov, G.V.Milovanivi\'{c} and Kh.M.Shadimetov \cite{HayShadMil14},  I.J.Schoenberg \cite{Schoen64},
L.L.Schumaker \cite{Schum81}, S.L.Sobolev \cite{Sob61a}, V.A.Vasilenko \cite{Vas83} and others.

The present work is also devoted to the variational method of construction of optimal interpolation formulas.

Assume, we are given the table of the values $\varphi(x_\beta)$, $\beta=0,1,...,N$ of a
function $\varphi$ at the points $x_{\beta}\in [0,1]$. It is required to approximate the function $\varphi$ by
another more simple function $P_\varphi$, i.e.
\begin{equation}\label{(1)}
\varphi (x) \cong P_\varphi  (x) = \sum\limits_{\beta = 0}^N
{C_\beta(x) \cdot \varphi (x_\beta )}
\end{equation}
which satisfies the following interpolation conditions
$$
\varphi(x_{\beta})=P_{\varphi}(x_{\beta}), \ \beta=0,1,...,N.
$$
Here $C_\beta (x)$ and $x_\beta$ ($\in [0,1]$) are
\emph{the coefficients} and \emph{the nodes} of the interpolation formula (\ref{(1)}), respectively.
By ${K_2}(P_2)$ we denote the class of all functions $\varphi$ defined on [0,1] which posses an absolutely continuous first derivative on [0,1] and whose second derivative is in $L_2(0,1)$. The class $K_2(P_2)$ under the pseudo-inner product
\begin{equation}\nonumber
\langle {\varphi,\psi} \rangle  =
\int\limits_0^1 \left(\varphi ''(x)+\omega^2\varphi(x)\right)\left(\psi_{\ell}''(x)+\omega^2\psi_{\ell}(x)\right)dx.
\end{equation}
is a Hilbert space if we identify functions that differ by $\sin(\omega x)$ and $\cos(\omega x)$, where $\omega\neq 0$.
Here we consider the norm
\begin{equation}\nonumber
\left\| {\varphi|K_2(P_2)} \right\| = \left\{
{\int\limits_0^1 {\left( \varphi''(x)+\omega^2\varphi(x)
\right)} ^2 \d x} \right\}^{1/2}.
\end{equation}
For the fixed $z\in [0,1]$ the error
\begin{equation}\label{(3)}
(\ell,\varphi)=\int\limits_{-\infty}^{\infty}\ell(x)\varphi(x)dx=\varphi(z)-P_{\varphi}(z)=\varphi(z)-\sum\limits_{\beta = 0}^N
{C_\beta(z) \cdot \varphi (x_\beta )}
\end{equation}
of the interpolation formula (\ref{(1)}) is a linear functional.\\
Here, in fixed $z\in [0,1]$,
$$
\ell(x,z)=\delta(x-z)-\sum_{\beta=0}^NC_\beta(z)\delta(x-x_\beta),
$$
and is the error functional of the
interpolation formula (\ref{(1)}) belongs to the space $K_2^*(P_2)$. Here $K_2^*(P_2)$ is the
conjugate space to the space $K_2(P_2)$, $\delta$ is  Dirac's delta-function.

By the Cauchy-Schwarz inequality the absolute value of the error (\ref{(3)}) is
estimated as follows
$$
|(\ell,\varphi)|\leq
\|\varphi|{K_2(P_2)}\|\cdot
\|\ell|{K_2^{*}(P_2)}\|,
$$
where
$$
\left\| {\ell|{K_2^{*}(P_2)} } \right\| = \mathop {\sup }\limits_{\varphi,\ \|\varphi \| \neq 0}
\frac{|(\ell,\varphi)|}{\|\varphi|{K_2(P_2)}\|}.
$$

Therefore, in order to estimate the error of the interpolation formula (\ref{(1)}) on
functions of the space $K_2(P_2)$ it is required to find the norm of the error functional $\ell$
in the conjugate space $K_2^{*}(P_2)$. That is we get the following problem.

\medskip

\textbf{Problem 1.1.} \textit{Find the norm of the error functional $\ell$ of the interpolation
formula (\ref{(1)}) in the space $K_2^{*}(P_2)$.}

\medskip

It is clear that the norm of the error functional $\ell$ depends on the coefficients $C_\beta(z)$
and the nodes $x_\beta$. The problem of minimization of the quantity $\|\ell\|$ by
coefficients $C_\beta(z)$ is the linear problem and by nodes $x_\beta$ is, in general, nonlinear
and complicated problem. We consider the problem of minimization of the
quantity $\|\ell\|$ by coefficients $C_\beta(z)$ when nodes $x_\beta$ are fixed.

The coefficients $\mathring{C}_{\beta}(z)$ (if there exist) satisfying the  equality
\begin{equation}\label{(4)}
\left\| \mathring{\ell}|{K_2^{*}(P_2)}\right\| = \mathop {\inf }\limits_{C_\beta(z) }  \left\|\ell|{K_2^{*}(P_2)} \right\|
\end{equation}
are called \textit{the optimal coefficients} and corresponding interpolation formula
$
\mathring{P}_{\varphi}(z)=\sum\limits_{\beta=0}^N\mathring{C}_{\beta}(z)\varphi(x_\beta)
$
is called \textit{the optimal interpolation formula} in the space $K_2^{*}(P_2)$.
Therefore, for construction of the interpolation formula we should solve the next problem.

\medskip

\textbf{Problem 1.2.} \textit{Find the coefficients $\mathring{C}_{\beta}(z)$ which satisfy equality (\ref{(4)}) when the
nodes $x_\beta$ are fixed.}

\medskip

It should be noted that Problems 1.1 and 1.2 were solved in \cite{BabHay18} in the Hilbert space $W_2^{(m,m-1)}$.
There the optimal interpolation formulas, which are exact for any polynomial of degree $m-2$ and for the function $\exp(-x)$, were obtained.

The rest of the paper is organized as follows. In Section 2, using the extremal function, the norm of the error functional is found.
Existence and uniqueness of the optimal interpolation formula of the form (\ref{(1)}) is discussed in Section 3. Section 4 is devoted to construction of the optimal interpolation formula. Finally, in Section 5 some numerical results are presented.

\section{The extremal function and the norm of the error functional $\ell$ }

Here we find explicit form of
the norm of the error functional $\ell$.

For finding the explicit form of the norm of the error functional
$\ell$ in the space $K_2(P_2)$ we use its extremal function which was introduced by Sobolev
\cite{Sob61a,Sob74}. The function  $\psi_{\ell}$ from
${K_2}(P_2)$ space is called \emph{the extremal
function} for the error functional $\ell$ if the following
equality is fulfilled
$$
\left( {\ell,\psi_{\ell}} \right) = \left\| {\ell\left|
{{K_2^*(P_2)}} \right.} \right\| \cdot \left\|
{\psi_{\ell}\left|{{K_2(P_2)}} \right.} \right\|.
$$
According to the Riesz theorem any linear continuous functional
$\ell$ in a Hilbert space is represented in the form of a
inner product. So, in our case, for any function $\varphi$ from ${K_2}(P_2)$
space we have
\begin{equation} \label{(5)}
\left( {\ell,\varphi} \right) = \langle{\psi _\ell,\varphi} \rangle.
\end{equation}
Here $\psi_\ell$ is the function from
${K_2}(P_2)$ is defined  uniquely by functional $\ell$ and is the extremal function.

It is easy to see from (\ref{(5)}) that the error functional $\ell$, defined on the space $K_2(P_2)$, satisfies the following equalities
\begin{eqnarray}
&&(\ell,\sin(\omega x))=0, \label{(6)}\\
&&(\ell,\cos(\omega x))=0.    \label{(7)}
\end{eqnarray}
The equalities (\ref{(6)}) and (\ref{(7)}) mean that our interpolation formula is exact for the functions $\sin(\omega x)$ and $\cos(\omega x)$.

The equation (\ref{(5)}) was solved in \cite{HayShadMil14} and for the extremal function $\psi_{\ell}$ was obtained the following
expression
\begin{equation}\label{(8)}
\psi_{\ell}(x)=(\ell*G_2)(x)+d_1\sin(\omega x)+d_2\cos(\omega x),
\end{equation}
where
\begin{equation}\label{(9)}
G_2(x)=\frac{\mathrm{sgn}x}{4\omega^3}[\sin(\omega x)-\omega x\cos(\omega x)],
\end{equation}
$*$ is the operation of convolution which for the functions $f$ and $g$ is defined as follows
\begin{equation}\label{(10)}
(f*g)(x)=\int\limits_{-\infty}^{\infty}f(x-y)g(y)\d y=\int\limits_{-\infty}^{\infty}f(y)g(x-y)\d y.
\end{equation}
Now we obtain the norm of the error functional
$\ell$. Since the space ${K_2}(P_2)$ is the Hilbert
space then by the Riesz theorem we have
\begin{equation}\label{(11)}
 \left( {\ell,\psi_\ell} \right) = \|\ell\| \cdot \| \psi _\ell\| =\| \ell\|^2 .
\end{equation}
Hence, using (\ref{(8)}) and (\ref{(9)}), taking
into account (\ref{(10)}) and (\ref{(11)}),  we get
{\small
\begin{eqnarray*}
\|\ell\|^2&=&
(\ell,\psi_{\ell})=\int\limits_{-\infty}^{\infty}\ell(x,z)\psi_{\ell}(x)\d x\\
&=&\int\limits_{-\infty}^{\infty}\left(\delta(x-z)-\sum\limits_{\beta=0}^NC_{\beta}(z)\delta(x-x_{\beta})\right)\\
&&\times \left(G_2(x-z)-\sum\limits_{\beta=0}^NC_{\beta}(z)G_2(x-x_{\beta})\right)\d x.
\end{eqnarray*}}
Hence, keeping in mind that $G_2(x)$, defined by (\ref{(9)}), is the even function, we
have
\begin{equation}
\|\ell\|^2=(-1)^m\left(\sum\limits_{\beta=0}^N\sum\limits_{\gamma=0}^N
C_{\beta}(z)C_{\gamma}(z)G_2(x_\beta-x_\gamma)-2\sum\limits_{\beta=0}^NC_{\beta}(z)G_2(z-x_{\beta})
\right).\label{(12)}
\end{equation}
Thus, Problem 1.1 is solved.

Further, we solve Problem 1.2.

\section{Existence and uniqueness of the optimal interpolation formula}

Assume that the nodes $x_\beta$  of the interpolation formula
(\ref{(1)}) are fixed. The error functional (\ref{(3)})
satisfies the conditions (\ref{(6)}) and (\ref{(7)}). The
norm of the error functional $\ell$ is  a multivariable
function with respect to the coefficients $C_\beta(z)$ $(\beta  =
\overline {0,N} )$. For finding the point of the conditional
minimum of the expression (\ref{(12)}) under the conditions
(\ref{(6)}) and (\ref{(7)}) we apply the Lagrange
method.

Consider the function
\begin{eqnarray*}
&&\Psi (C_0(z),C_1(z),...,C_N(z),d_1,d_2)\\
&&\qquad=\left\| \ell \right\|^2  - 2\Bigg(d_1(z)(\ell,\sin(\omega x))+d_2(z)(\ell,\cos(\omega x))\Bigg).
\end{eqnarray*}
Equating to 0 the partial derivatives of the function $\Psi$
by $C_\beta(z)$ $(\beta =\overline{0,N})$, $d_1(z)$ and $d_2(z)$, we get the following
system of $N+3$ linear equations of $N+3$ unknowns
\begin{eqnarray}
&&\sum\limits_{\gamma=0}^N C_\gamma(z) G_2(x_\beta   - x_\gamma  )
+ d_1(z)\sin(\omega x_\beta)+d_2(z)\cos(\omega x_\beta) = G_2(z-x_{\beta}),\nonumber\\
&&\beta  =0,1,...,N,\    \label{(13)} \\
&& \sum\limits_{\gamma=0}^N C_\gamma(z)\sin(\omega x_\gamma)=\sin(\omega z),\ \label{(14)}\\
&&\sum\limits_{\gamma=0}^NC_\gamma(z)\cos(\omega x_\gamma)=\cos(\omega z),\label{(15)}
\end{eqnarray}
where $G_2(x)$ is defined by equality (\ref{(9)}).

The system (\ref{(13)})-(\ref{(15)}) has a unique solution
and this solution gives the minimum to $\left\|\ell\right\|^2
$ under the conditions (\ref{(14)}) and (\ref{(15)}).

The uniqueness of the solution of the system
(\ref{(13)})--(\ref{(15)}) is proved as the uniqueness of the solution of the system (24)--(26) of the work
\cite{SobVas}.

Therefore, in fixed values of the nodes $x_\beta$ the square of
the norm of the error functional $\ell$, being quadratic
function of the coefficients $C_\beta(z)$, has a unique minimum in
some concrete value $C_\beta(z)= \mathring{C}_{\beta }(z)$.

\begin{rem}\label{Rem(1)}
It should be noted that by integrating both sides of the system (\ref{(13)})-(\ref{(15)}) by $z$ from 0 to 1
we get the system (3.1)-(3.1) of the work \cite{HayShadMil14}. This means that by integrating the optimal interpolation formula (\ref{(1)}) in  the space $K_2(P_2)$ we get the optimal quadrature formula of the form (\ref{(1)}) in the same space (see \cite{HayShadMil14}).
\end{rem}

\begin{rem} \label{Rem(2)}
It is clear from the system (\ref{(13)})-(\ref{(15)}) that for the optimal coefficients the following are true
\begin{eqnarray*}
\mathring{C}_{\beta}(h\gamma)&=&\left\{
\begin{array}{ll}
1, & \gamma =\beta,\\
0,& \gamma\neq \beta,
\end{array}
\right. \gamma=0,1,...,N,\ \ \beta=0,1,...,N.
\end{eqnarray*}
\end{rem}

Below for convenience the optimal coefficients $ \mathring{C}_{\beta}(z)$ we remain as $C_\beta(z)$.
\section{The algorithm for computation of coefficients of the optimal interpolation formula}

In the present section we give the algorithm for solution of the
system (\ref{(13)})-(\ref{(15)}). Below mainly is used the concept of discrete
argument functions and operations on them. The theory of discrete
argument functions is given, for instance, in \cite{Sob74,SobVas}. For
completeness we give some definitions about functions of discrete
argument.

Assume that the nodes $x_\beta$ are equal spaced, i.e. $x_\beta=
h\beta,$ $h = {1 \over N}$, $N = 1,2,...$.

\begin{defn}\label{def1}
  The function $\varphi (h\beta )$ is a
\emph{discrete argument function} (or \emph{discrete function}) if it is given on some set of
integer values of $\beta$.
\end{defn}

\begin{defn}\label{def2}
\emph{The inner product} of two discrete
functions $\varphi(h\beta )$ and $\psi (h\beta )$ is given by
$$
\left[ {\varphi(h\beta),\psi(h\beta) } \right] =
\sum\limits_{\beta  =  - \infty }^\infty  {\varphi (h\beta ) \cdot
\psi (h\beta )},
$$
if the series on the right hand side converges absolutely.
\end{defn}

\begin{defn}\label{def3}
\textit{The convolution} of two functions
$\varphi(h\beta )$ and $\psi (h\beta )$ is the inner product
$$
\varphi (h\beta )*\psi (h\beta ) = \left[ {\varphi (h\gamma ),\psi
(h\beta  - h\gamma )} \right] = \sum\limits_{\gamma  =  - \infty
}^\infty  {\varphi (h\gamma ) \cdot \psi (h\beta  - h\gamma )}.
$$
\end{defn}

Now we turn on to our problem.

Suppose that $C_\beta(z)=0$  when $\beta  < 0$ and $\beta  > N$.
Thus we have the following problem.

\medskip

\textbf{Problem 4.1} {\it
Find the discrete functions $C_\beta(z)$, $\beta=0,1,...,N$,  $d_1(z)$ and $d_2(z)$ which
satisfy the system (\ref{(14)})-(\ref{(15)}).}

\medskip

Further we investigate Problem 4.1 which is equivalent to Problem 1.2.
Instead of $C_\beta(z)$ we introduce the following functions
\begin{eqnarray}
v_2(h\beta ) &=& G_2(h\beta )*C_\beta(z),\nonumber\\
u_2(h\beta)&=&v_2\left({h\beta}\right)+ d_1(z)\sin(\omega h\beta)+d_2\cos(\omega h\beta).\label{(17)}
\end{eqnarray}

Now we should express the coefficients
$C_\beta(z)$ by the function $u(h\beta)$. For this we use the operator $D_2(h\beta)$ which satisfies the equality
\begin{equation}
D_2(h\beta )*G_2(h\beta ) =\delta(h\beta) ,\label{(18)}
\end{equation}
where $\delta(h\beta)$ is equal to 0 when $\beta \ne 0$ and is
equal to 1 when $\beta  = 0$, i.e. $\delta(h\beta)$ is the
discrete delta-function.

In \cite{Hay04a} the
operator $D_2(h\beta )$  which satisfies
equation (\ref{(18)}) is constructed and its some properties
are studied.

The following theorems are proved in \cite{Hay04a}.

\begin{thm}\label{THM4.1}
The discrete analogue of the differential operator ${{d^{4}}\over {d x^{4} }}+2\omega^2{{d^{2}} \over {d x^{2} }}+\omega^4$ satisfying the
equation (\ref{(18)}) has the form
\begin{equation}
D_2(h\beta)=p\left\{
\begin{array}{ll}
A\lambda^{|\beta|-1}, &|\beta|\geq 2,\\
1+A,& |\beta|=1,\\
C+\frac{A}{\lambda}, & \beta=0.
\end{array}
\right. \label{(19)}
\end{equation}
where
\begin{equation}
\displaystyle
\begin{array}{l}
p=\frac{2\omega^3}{\sin(\omega h)-\omega h\cos(\omega h)},\nonumber\\
C=\frac{2h\omega\cos(2h\omega)-\sin(2h\omega)}{\sin(h\omega)-h\omega\cos(h\omega)},\nonumber\\
A=\frac{(2h\omega)^2\sin^4(h\omega)\cdot\lambda^2}{(\lambda^2-1)(\sin(h\omega)-h\omega\cos(h\omega))^2},\nonumber\\
\lambda=\frac{2h\omega-\sin(2h\omega)-2\sin(h\omega)\cdot\sqrt{(h\omega)^2-\sin^2(h\omega)}}{2(h\omega\cos(h\omega)-\sin(h\omega))},\  |\lambda|<1.
\end{array}
\end{equation}

\end{thm}
\begin{thm}\label{THM4.2}
The discrete analogue $D_2(h\beta)$ of the
differential operator ${{d^{4}} \over {d x^{4} }}+2\omega^2{{d^{2}}\over {d x^{2} }}+\omega^4$ satisfies the following equalities

1) $D_2(h\beta)*\sin(h\omega\beta)=0,$

2) $D_2(h\beta)*\cos(h\omega\beta)=0,$

3) $D_2(h\beta)*(h\omega\beta)\cos(h\omega\beta)=0,$

4) $D_2(h\beta)*(h\omega\beta)\sin(h\omega\beta)=0.$\\
\end{thm}

Then taking into account (\ref{(18)}),  (\ref{(19)}), using
Theorems \ref{THM4.1} and \ref{THM4.2}, for optimal coefficients we
have
\begin{equation}
C_\beta(z)=D_2(h\beta )*u_2(h\beta ).\label{(28)}
\end{equation}

Thus if we find the function $u_2(h\beta )$ then the optimal
coefficients $C_{\beta}(z)$ will be found from equality
(\ref{(28)}).

In order to calculate the convolution (\ref{(28)}) it is required to
find the representation of the function $u_2(h\beta )$ at all
integer values of $\beta$. From equality (\ref{(17)}) we get that
$u_2(h\beta ) = G_2(z-h\beta )$ when $h\beta \in [0,1]$. Now we find the representation of the function $u_2(h\beta )$ when
$\beta < 0$ and $\beta>N$.

Since $C_\beta(z)= 0$ when $h\beta \notin [0,1]$ then
$$
C_\beta(z)   = D_2(h\beta )*u_2(h\beta ) =
0,\,\,\,\, h\beta \notin [0,1].
$$

Now we calculate the convolution $v_2(h\beta ) =G_2(h\beta)*C_\beta(z)$ when $h\beta \notin [0,1]$.

Suppose $\beta  \leq 0$ and $\beta\geq N$ then taking into account equalities
(\ref{(9)}), (\ref{(14)}) and (\ref{(15)}), we have
\begin{equation}
u_2(h\beta)=\left\{
\begin{array}{ll}
d_1^-(z)\sin(\omega h\beta)+d_2^-(z)\cos(\omega h \beta)+\frac{h\beta}{4\omega^2}\cos(\omega (h\beta-z)),& \beta\leq 0,\\[2mm]
G_2(z-h\beta),& 0\leq \beta\leq N,\\[2mm]
d_1^+(z)\sin(\omega h\beta)+d_2^+(z)\cos(\omega h\beta)-\frac{h\beta}{4\omega^2}\cos(\omega(h\beta-z)),& \beta\geq N,\\
\end{array}
\right. \label{(29)}
\end{equation}
where $d_1^-(z),\ d_1^+(z),\ d_2^-(z)$ and $d_2^+(z)$ are unknowns.

From (\ref{(29)}) when $\beta=0$ and $\beta=N$ we get
\begin{eqnarray}
&&d_2^-(z)=G_2(z),  \ \label{(30)}\\
&&d_2^+(z)=\frac{1}{\cos(\omega)}\Bigg[G_2(z-1)+\frac{\cos(\omega(1-z))}{4\omega^2}-d_1^+(z)\sin(\omega)\Bigg].\ \label{(31)}
\end{eqnarray}
Thus, putting (\ref{(30)}) and (\ref{(31)}) to (\ref{(29)}) we have the following explicit form of the function $u_2(h\beta)$:
\begin{equation}
u_2(h\beta)=\left\{
\begin{array}{ll}
d_1^-(z)\sin(\omega h\beta)+G_2(z)\cos(\omega h \beta)+\frac{h\beta}{4\omega^2}\cos(\omega (h\beta-z)),& \beta\leq 0,\\[2mm]
G_2(z-h\beta),& 0\leq \beta\leq N,\\[2mm]
d_1^+(z)\frac{\sin(\omega (h\beta-1))}{\cos(\omega)}+\frac{\cos(\omega h\beta)}{\cos(\omega)}\Bigg[G_2(z-1)+\frac{\cos(\omega(1-z) )}{4\omega^2}\Bigg]\\
-\frac{h\beta}{4\omega^2}\cos(\omega (h\beta-z)),& \beta\geq N.\\
\end{array}
\right. \label{(32)}
\end{equation}
In the last expression of the function $u_2(h\beta)$ we have only two unknowns $d_1^-(z)$ and $d_1^+(z)$.

Hence using (\ref{(19)}) and (\ref{(32)}) we get the following problem.\\

\medskip

\textbf{Problem 4.2.} {\it
Find the solution of the equation
\begin{equation}
D_2(h\beta )*u_2(h\beta) = 0,\,\, h\beta  \notin[0,1],\nonumber
\end{equation}
having the form (\ref{(32)}).
Here $d_1^{-}(z)$ and $d_1^{+}(z)$ are unknowns.}

\medskip

Unknowns $d_1^{-}(z)$ and $d_1^{+}(z)$ we find from the equation
\begin{equation}\label{(34)}
D_2(h\beta )*u_2(h\beta)=0
\end{equation}
when $\beta=-1$ and $\beta=N+1$. From the last equation

From the system (\ref{(34)}) in the case $\beta=-1$ and after some simplifications we have the following system of equations
\begin{equation}
A_{11}d_1^-(z)+A_{12}d_1^+(z)=S_1 \ \label{(35)}
\end{equation}
where
\begin{equation}
\begin{array}{ll}
A_{11}=-\sin(2\omega h)-C\sin(\omega h)-\frac{A\sin(\omega h)}{\lambda(\lambda^2-2\lambda\cos(\omega h)+1)},\\ \nonumber
A_{12}=\frac{A\lambda^{N+1}\sin(\omega h)}{\cos(\omega)\left(\lambda^2-2\lambda\cos(\omega h)+1\right)},\\
S_1=\frac{1}{4\omega^2}\Bigg(Ch\cos(\omega(h+z)+2h\cos(\omega(2h+z))-\frac{A\lambda^{N}\cos(\omega(1-z))}{\cos(\omega)})K_1\\
\ \ +\frac{Ah}{\lambda^2}K_2+Ah\lambda^NK_3\Bigg)-AK_4-\frac{A\lambda^NG_2(z-1)}{\cos(\omega)}K_1\\
\ \ -G_2(z)\Bigg(C\cos(\omega h)+1+\cos(2\omega h)+\frac{A}{\lambda^2}K_5\Bigg),\\
K_1=\frac{\lambda(\cos(\omega(1+h))-\lambda\cos(\omega))}{\lambda^2-2\lambda\cos(\omega h)+1},\\
K_2=\frac{\lambda(\lambda^2\cos(\omega(h-z))-2\lambda\cos(\omega z)+\cos(\omega(h+z)))}{(\lambda^2-2\lambda\cos(\omega h)+1)^2},\\
\end{array}
\end{equation}
\begin{equation}
\begin{array}{ll}
K_3=\sum\limits_{\gamma = 1}^\infty\lambda^\gamma(N+\gamma)\cos(\omega (h(N+\gamma)-z)),\\ \nonumber
K_4=\sum\limits_{\gamma = 0}^N\lambda^\gamma G_2(z-h\gamma),\\
K_5=\frac{\lambda(\cos(\omega h)-\lambda)}{\lambda^2-2\lambda\cos(\omega h)+1}.\\
\end{array}
\end{equation}

Now, from (\ref{(34)}) in the case $\beta = N + 1$ doing some calculations we get the next equation
\begin{equation}
A_{21}d_1^-(z)+A_{22}d_1^+(z)=S_2,  \ \label{(36)}
\end{equation}
where
\begin{equation}
\begin{array}{ll}
A_{21}=-\frac{A\lambda^{N+1}\sin(\omega h)}{\lambda^2-2\lambda\cos(\omega h)+1},\\ \nonumber
A_{22}=\frac{1}{\cos(\omega)}\Bigg[C\sin(\omega h)+\sin(2\omega h)+\frac{A\sin(\omega h)}{\lambda(\lambda^2-2\lambda\cos(\omega h)+1)}\Bigg],\\
S_2=\frac{1}{4\omega^2}\Bigg(A\lambda^NhK_2+\cos(\omega(1-z))+Ch(N+1)\cos(\omega(h(N+1)-z))\\
\ \ +h(N+2)\cos(\omega(h(N+2)-z))+\frac{AhK_3}{\lambda^2}-\frac{\cos(w(1-z))K_6}{\cos(\omega)}\Bigg)\\
\ \ -A\lambda^NK_7-G_2(z)A\lambda^NK_5-\frac{G_2(z-1)}{\cos(\omega)}K_6,\\
K_6=\cos(\omega)+C\cos(\omega h(N+1))+\cos(\omega h(N+2))+\frac{AK_1}{\lambda^2},\\
K_7=\sum\limits_{\gamma = 0}^N\lambda^{-\gamma} G_2(z-h\gamma),\\
K_8=\frac{\lambda(\sin(\omega (h+1))-\lambda\sin(\omega))}{\lambda^2-2\lambda\cos(\omega h)+1}.\\ \nonumber
\end{array}
\end{equation}
Then solving the system  (\ref{(35)}), (\ref{(36)}) of equations we get $d_1^-(z)$ and $d_1^+(z)$.
Finally, from (\ref{(28)}) for $\beta=0,1,...,N$ we get the explicit formulas for optimal coefficients as claimed in the following theorem.

\begin{thm}\label{THM4.3}
Coefficients of the optimal interpolation formula (\ref{(1)})
with equal spaced nodes in the space $K_2(P_2)$ have the
following form
\begin{equation}
\begin{array}{ll}
C_{0}(z)=&
p\Bigg[-d_1^-(z)\sin(\omega h)+d_2^-(z)\cos(\omega h)-\frac{1}{4\omega^2}h\cos(\omega(h+z))+C G_2(z)\\\nonumber
&+G_2(z-h)+\frac{A}{\lambda}\Bigg[\sum\limits_{\gamma=0}^N\lambda^\gamma G_2(z-h\gamma)+M_1+\lambda^N N_1\Bigg]\Bigg],\\
C_{\beta}(z)=&
p\Bigg[G_2(z-h(\beta-1))+C G_2(z-h\beta)+G_2(z-h(\beta+1))\\ \nonumber
\end{array}
\end{equation}
\begin{equation}
\begin{array}{ll}
&+\frac{A}{\lambda}\Bigg[\sum\limits_{\gamma=0}^N\lambda^{|\beta-\gamma|} G_2(z-h\gamma)+\lambda^\beta M_1+\lambda^{N-\beta} N_1\Bigg]\Bigg],  \beta=1,2,..,N-1.\\ \nonumber
\end{array}
\end{equation}
\begin{equation}
\begin{array}{ll}\label{(37)}
C_{N}(z)=&p\Bigg[d_1^+(z)\sin(\omega h(N+1))+d_2^+(z)\cos(\omega h(N+1))\\ \nonumber
&-\frac{1}{4\omega^2}h(N+1)\cos(\omega(h(N+1)-z))+C G_2(z-1)\\
&+G_2(z-h(N-1))+\frac{A}{\lambda}\Bigg[\sum\limits_{\gamma=0}^N\lambda^{N-\gamma} G_2(z-h\gamma)+\lambda^N M_1+ N_1\Bigg]\Bigg],\\
\end{array}
\end{equation}
here
\begin{equation}
\begin{array}{ll}
M_1=-d_1^-(z)\frac{\lambda\sin(\omega h)}{\lambda^2-2\lambda\cos(\omega h)+1}+d_2^-(z)K_5-\frac{h}{4\omega^2}K_2, \\ \nonumber
N_1=d_1^+(z)K_8+d_2^+(z)K_1-\frac{h}{4\omega^2} K_3,\\
d_1^-(z)=\frac{S_1\cdot A_{22} - S_2 \cdot A_{12}}{A_{11} \cdot A_{22} - A_{21} \cdot A_{12} },\\
d_1^+(z)=\frac{S_2\cdot A_{11} - S_1 \cdot A_{21}}{A_{11} \cdot A_{22} - A_{21} \cdot A_{12} },\\
d_2^-(z)=G_2(z),\\
d_2^+(z)=\frac{1}{\cos(\omega)}\Bigg[G_2(z-1)+\frac{\cos(\omega(1-z))}{4\omega^2}-d_1^+(z)\sin(\omega)\Bigg].
\end{array}
\end{equation}
\end{thm}

\section{Numerical results}
In this section we give some numerical results.

First, when $N=5$ using Theorem \ref{THM4.3}, we get the graphs of the coefficients of the optimal interpolation formulas
$$
\varphi(z)\cong \mathring{P}_{\varphi}(z)=\sum\limits_{\beta=0}^5\mathring{C}_{\beta}(z)\varphi(h\beta),\ z\in [0,1].
$$ They are presented in Fig 1.
These graphical results confirm Remark \ref{Rem(2)} for the case $N=5$ , i.e. for the optimal coefficients the following hold
$$
\mathring{C}_{\beta}(h\gamma)=\delta_{\beta\gamma},\ \ \beta,\gamma=0,1,...,5,
$$
where $\delta_{\beta\gamma}$ is the Kronecker symbol.

Now, in numerical examples, we interpolate the functions
$$
\varphi_1(z)=z^2,\ \varphi_2(z)=\E^{z}\mbox{ and } \varphi_3(z)=\sin z
$$
by optimal interpolation formulas of the form (\ref{(1)}) in the cases $N=5,\ 10$, using Theorem \ref{THM4.3}. For the functions
$\varphi_i$, $i=1,2,3$ the graphs of absolute errors $|\varphi_i(z)-\mathring{P}_{\varphi_i}(z)|$, $i=1,2,3$, are given in Fig 2, Fig 3, Fig 4. In these Figures one can see that by increasing value of $N$ absolute errors between optimal interpolation formulas and given functions are decreasing.
 \vspace{-50pt}
\begin{center}
  \includegraphics[width=0.3\textwidth]{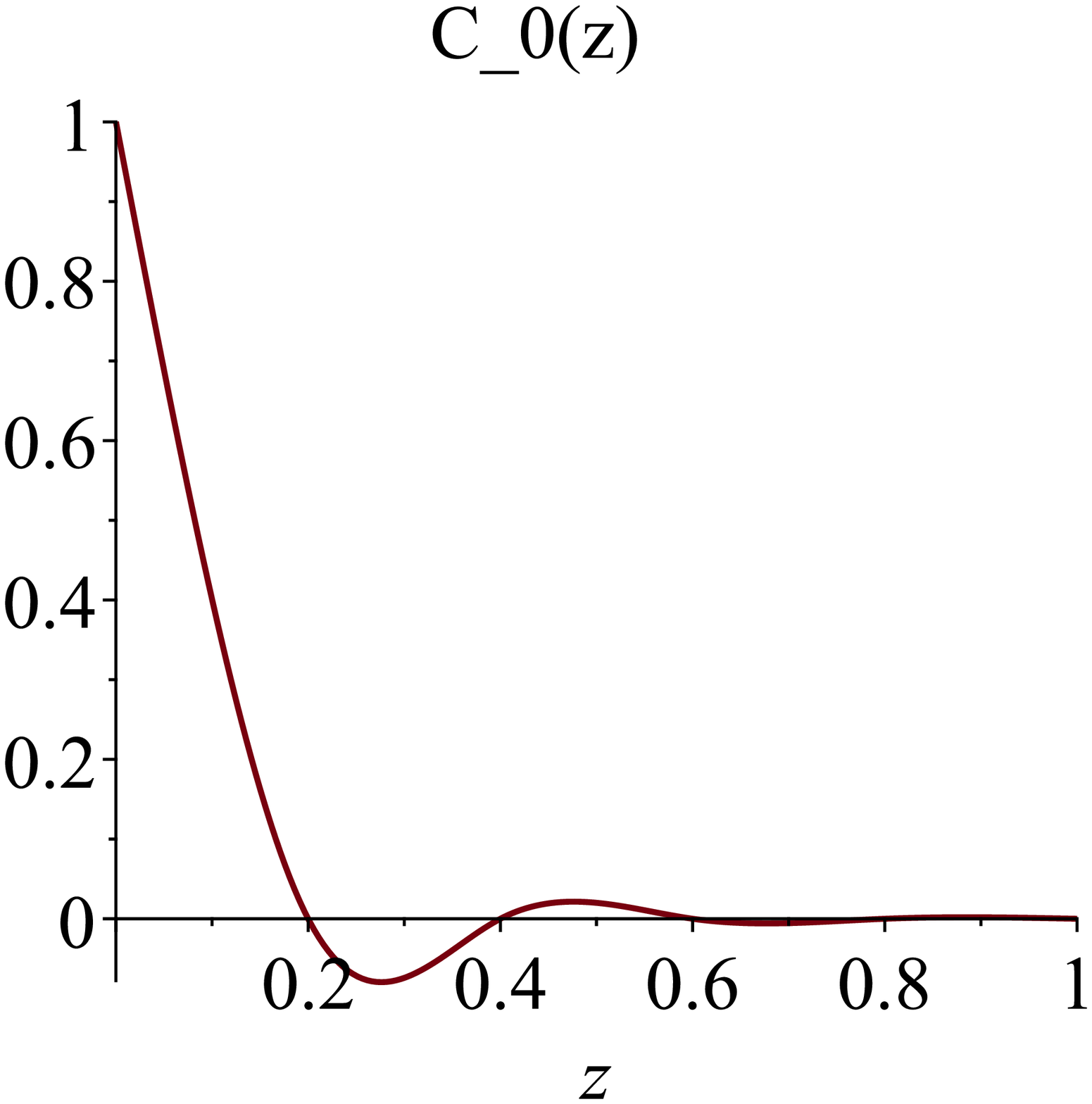}
  \includegraphics[width=0.3\textwidth]{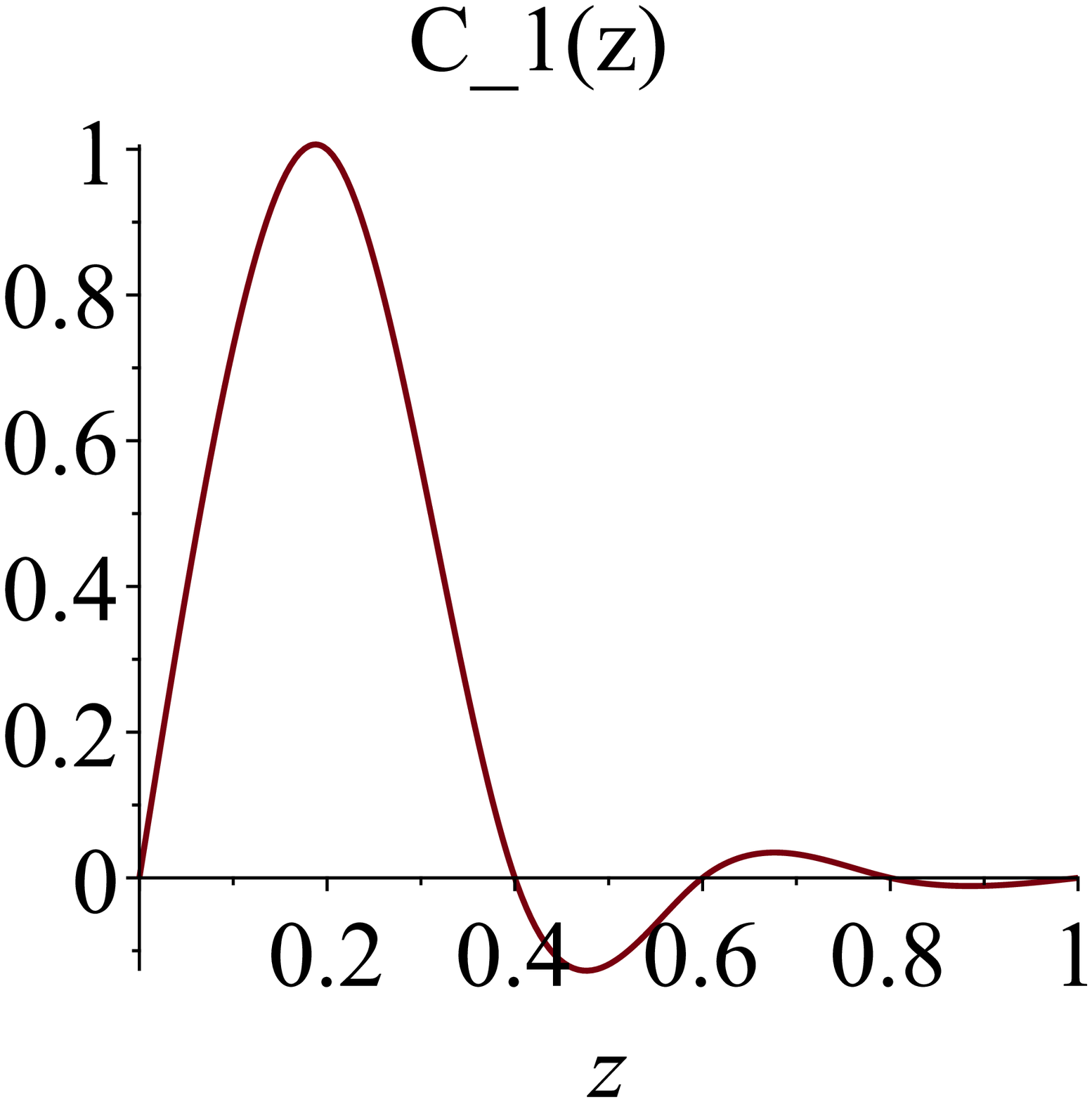}
  \includegraphics[width=0.3\textwidth]{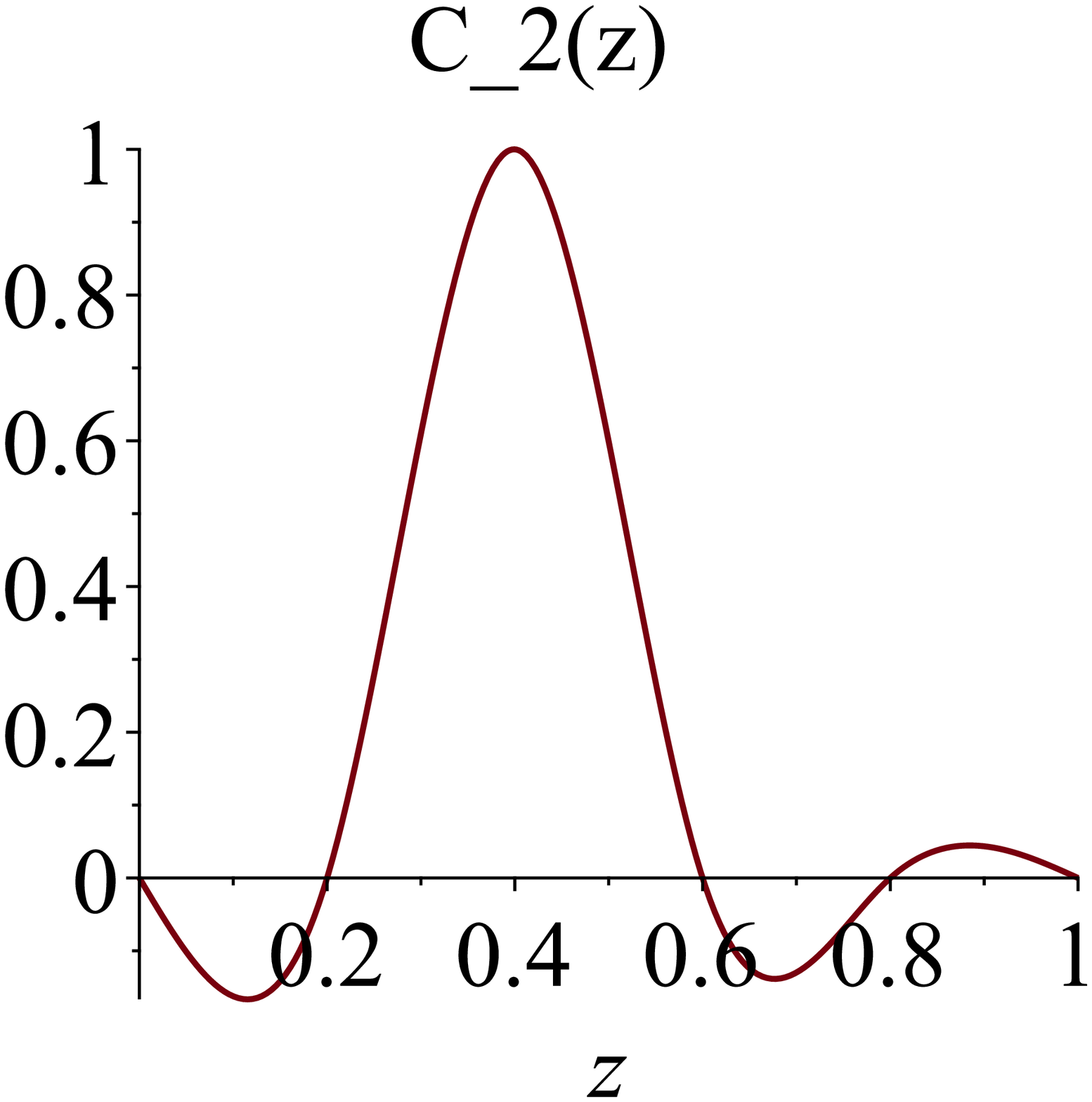}\\
  \vspace{-50pt}
  \includegraphics[width=0.3\textwidth]{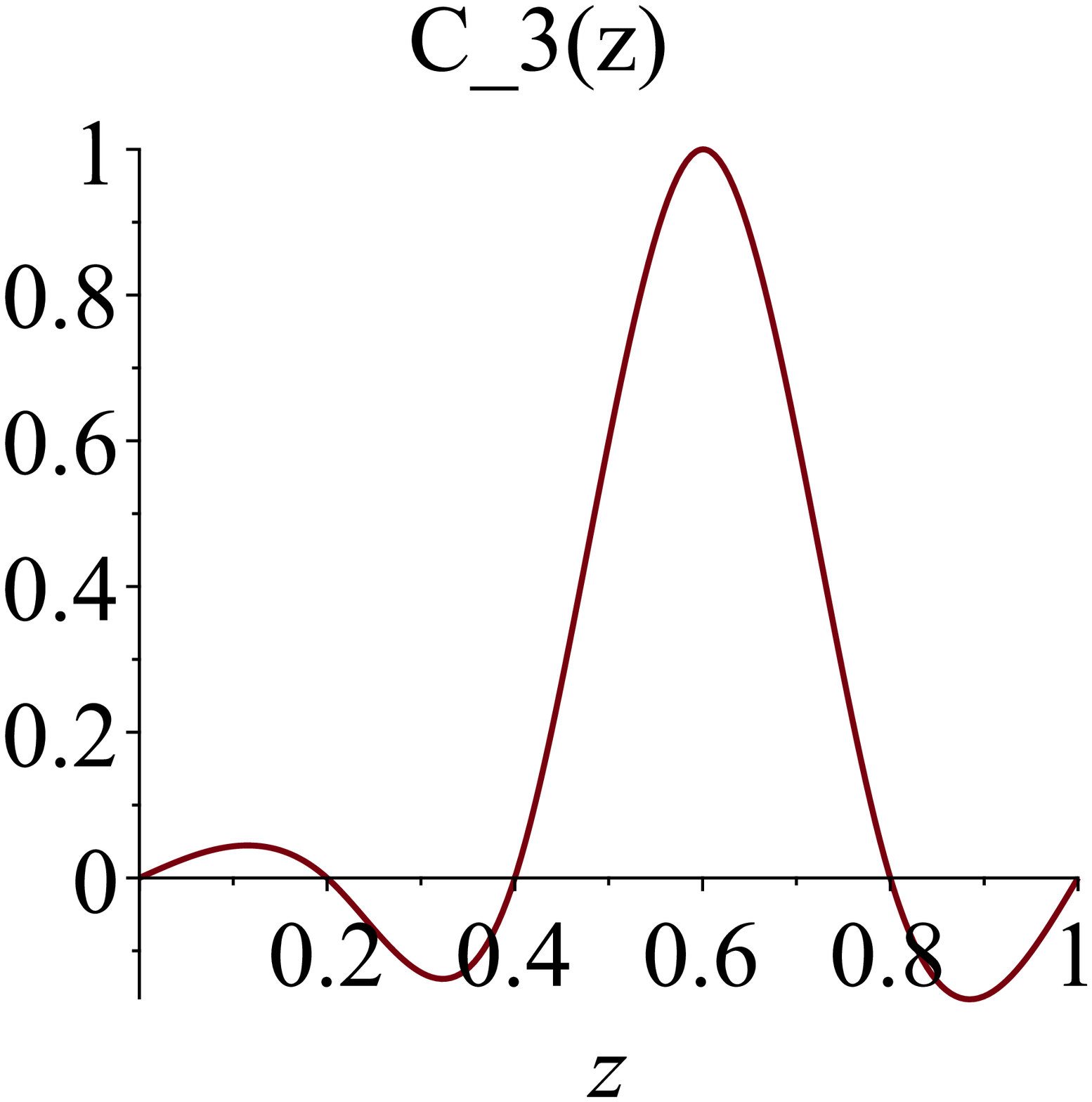}
  \includegraphics[width=0.3\textwidth]{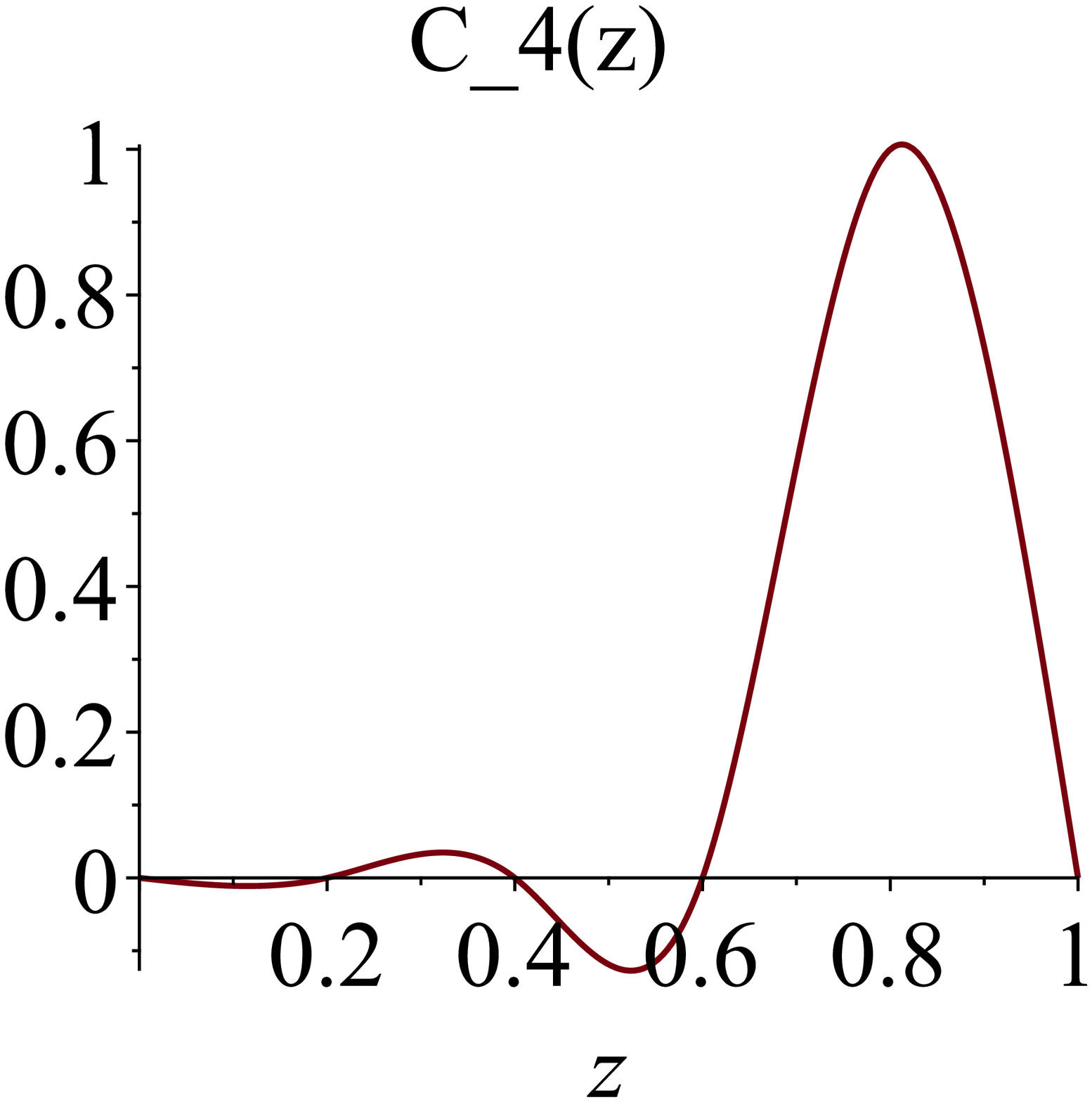}
  \includegraphics[width=0.3\textwidth]{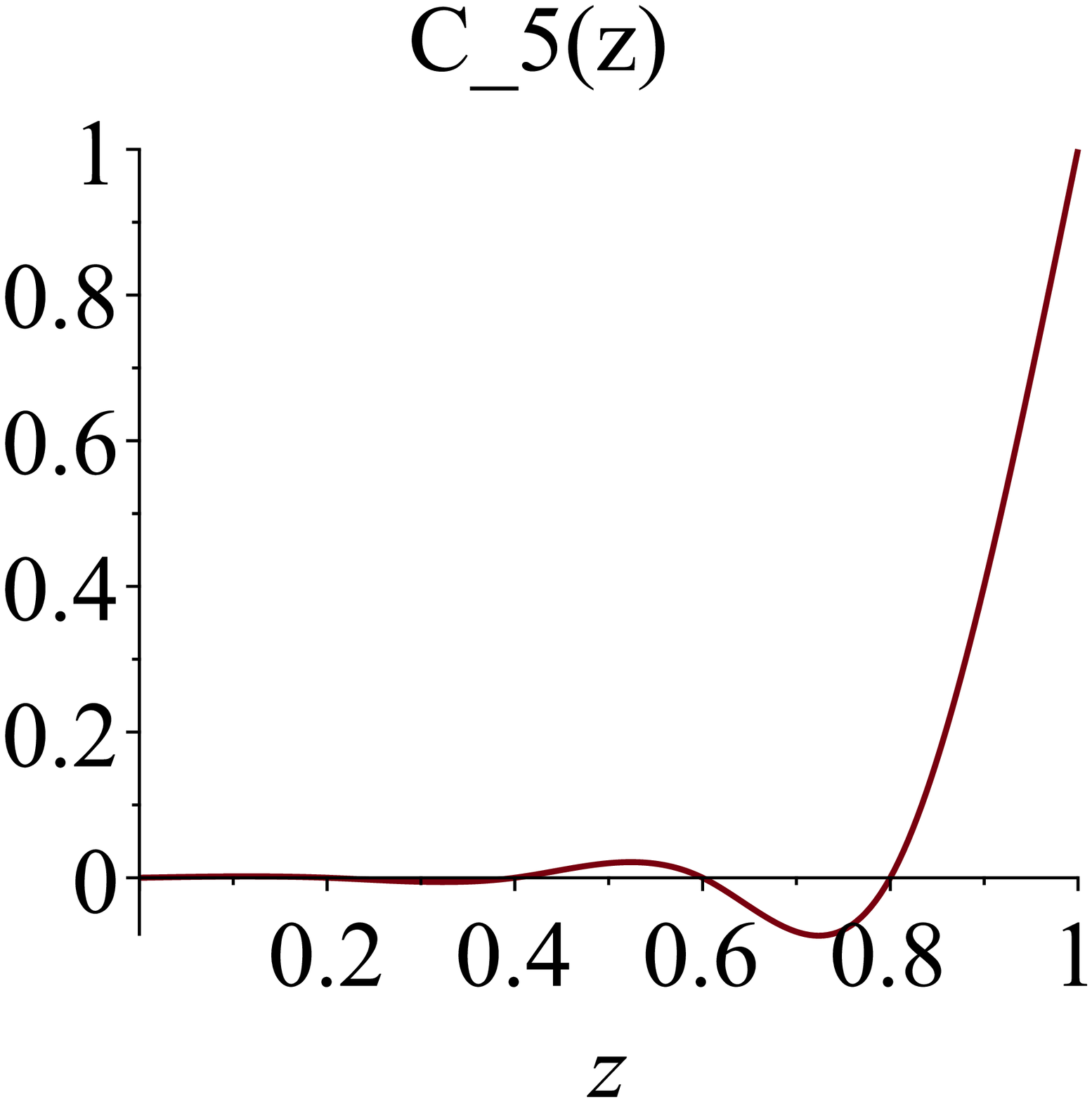}\\
   \vspace{-50pt}
Fig. 1. {\small Graphs of coefficients of the optimal interpolation formulas (\ref{(1)}) in the case  $N=5$.}
\end{center}
\begin{center}
  \includegraphics[width=0.45\textwidth]{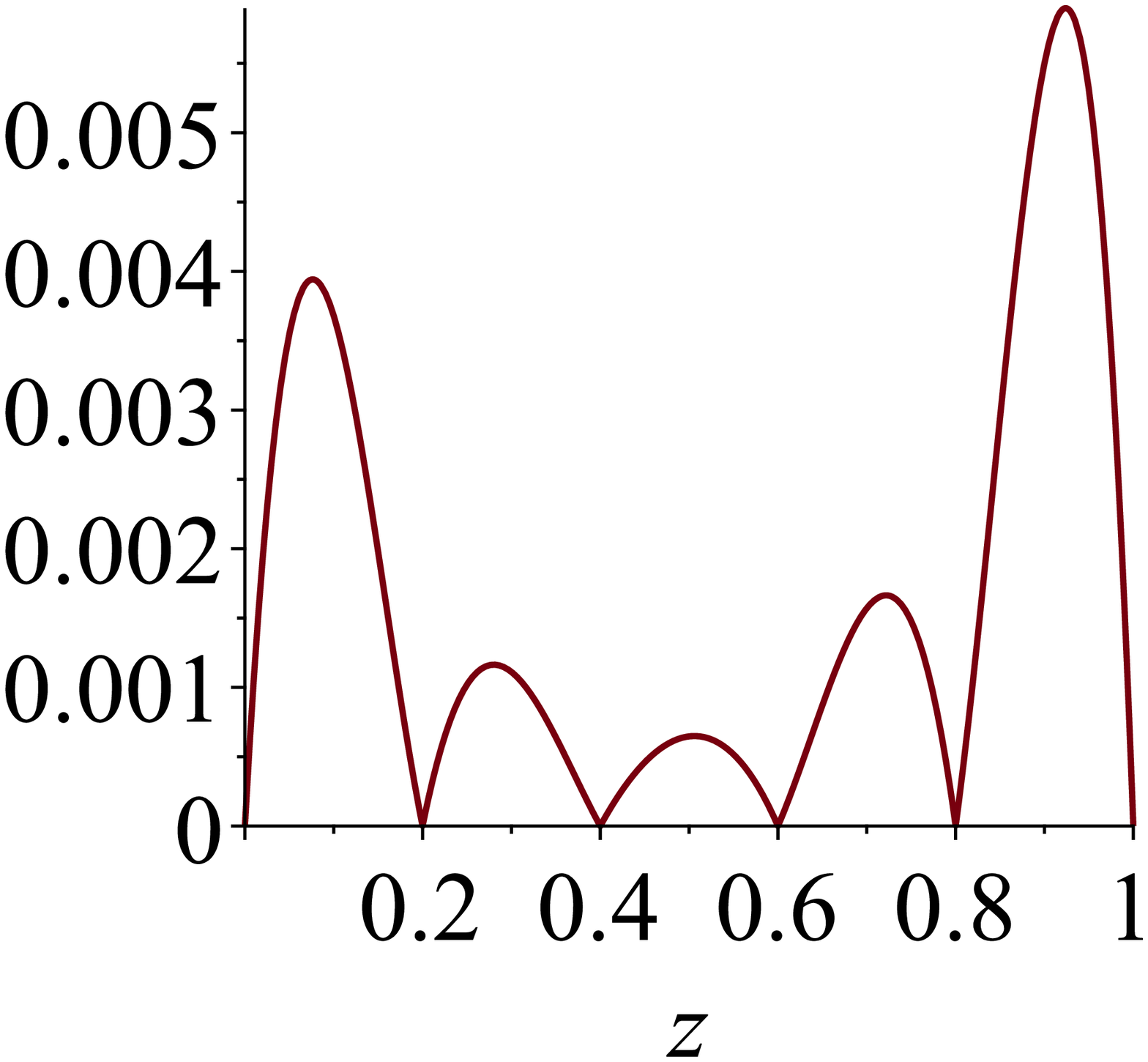}
  \includegraphics[width=0.45\textwidth]{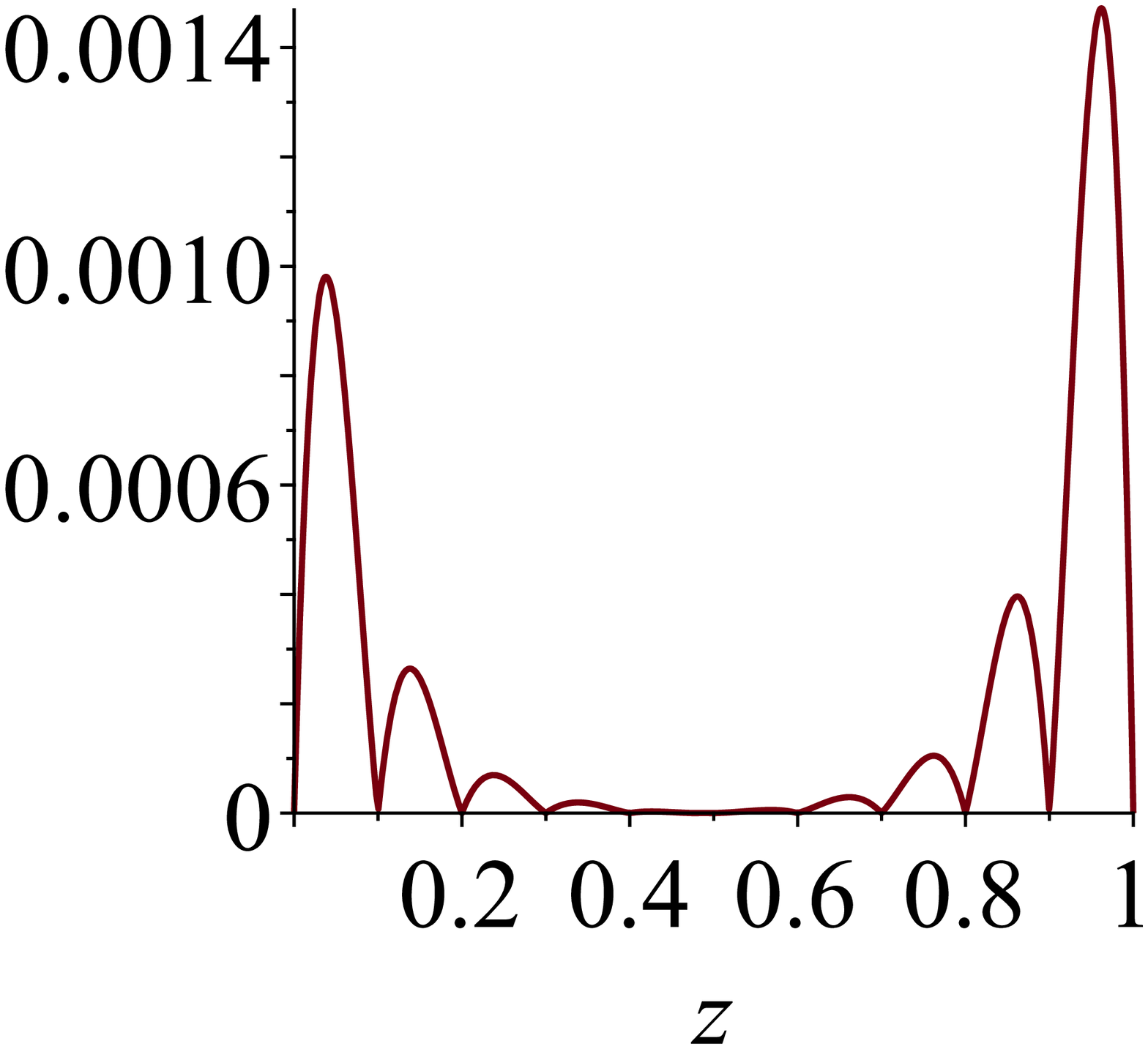}\\
   \vspace{-60pt}
Fig. 2. { Graphs of absolute errors for  $N=5$ and $N=10$:  $|z^2-P_{z^2}(z)|$.}
\end{center}
 \vspace{-40pt}
 \begin{center}
\includegraphics[width=0.45\textwidth]{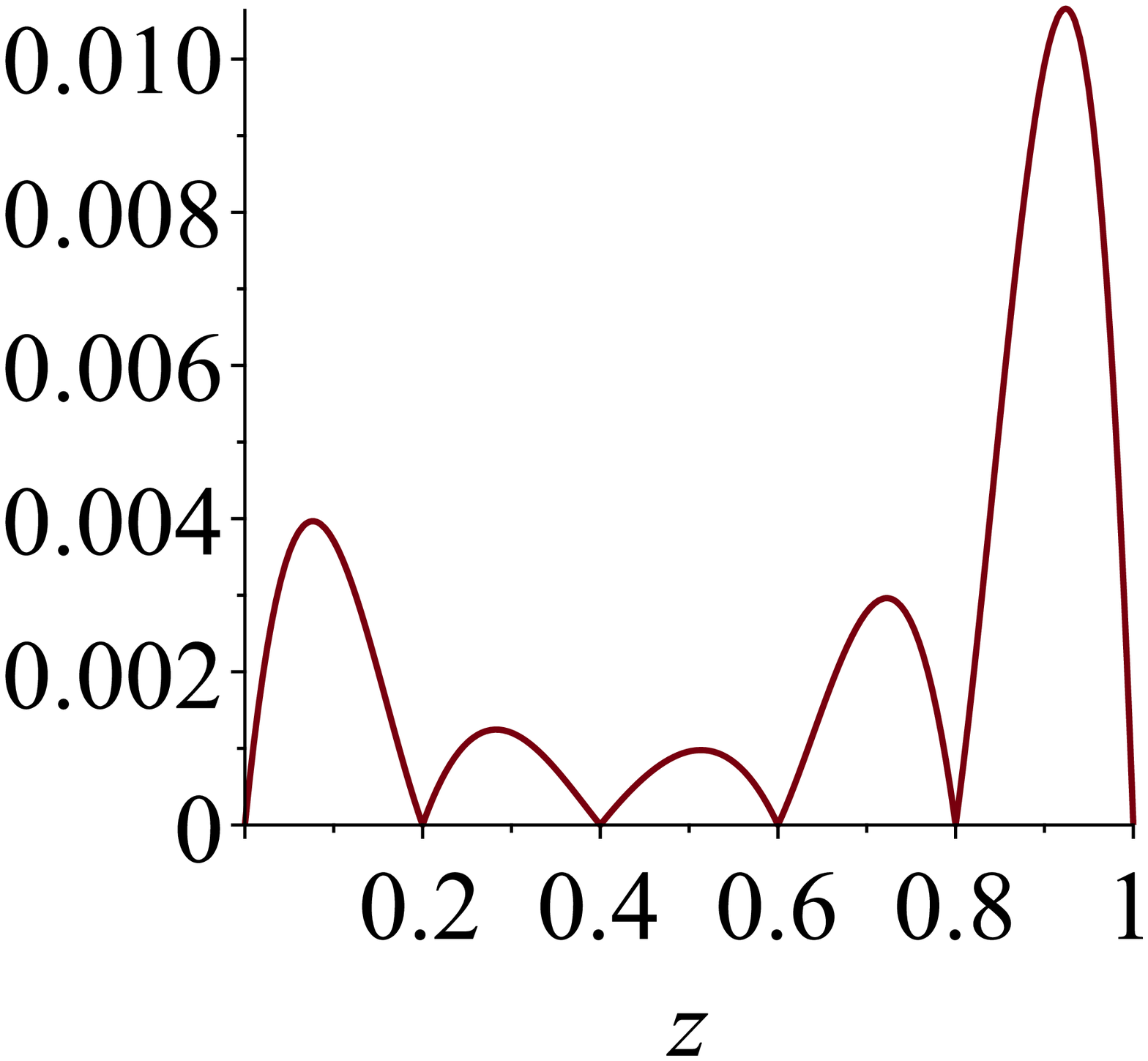}
  \includegraphics[width=0.45\textwidth]{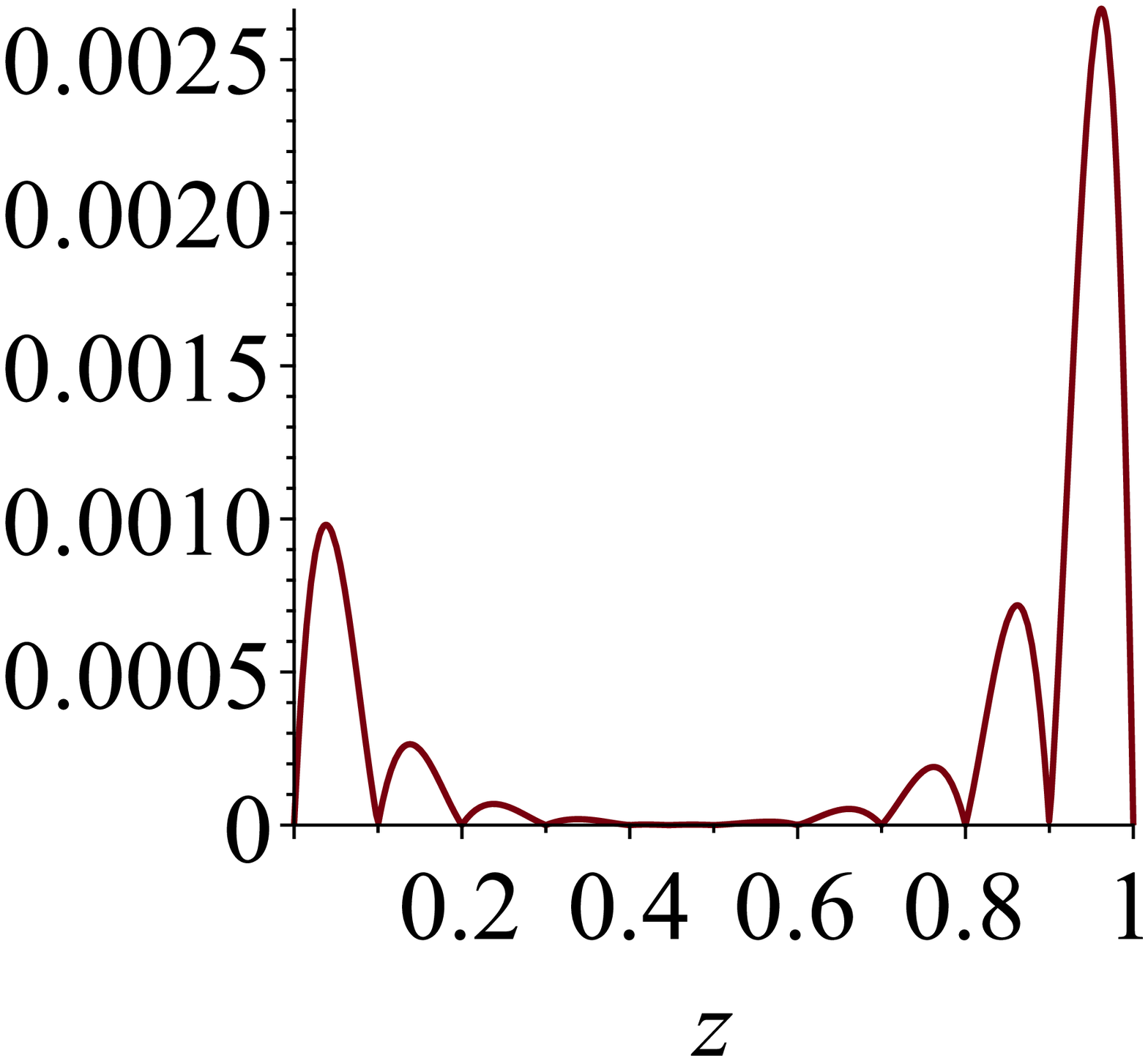}\\
   \vspace{-60pt}
Fig. 3. {\small Graphs of absolute errors for  $N=5$ and $N=10$:  $|\exp(z)-P_{\exp(z)}(z)|$.}
\end{center}
 \vspace{-30pt}
 \begin{center}
  \includegraphics[width=0.45\textwidth]{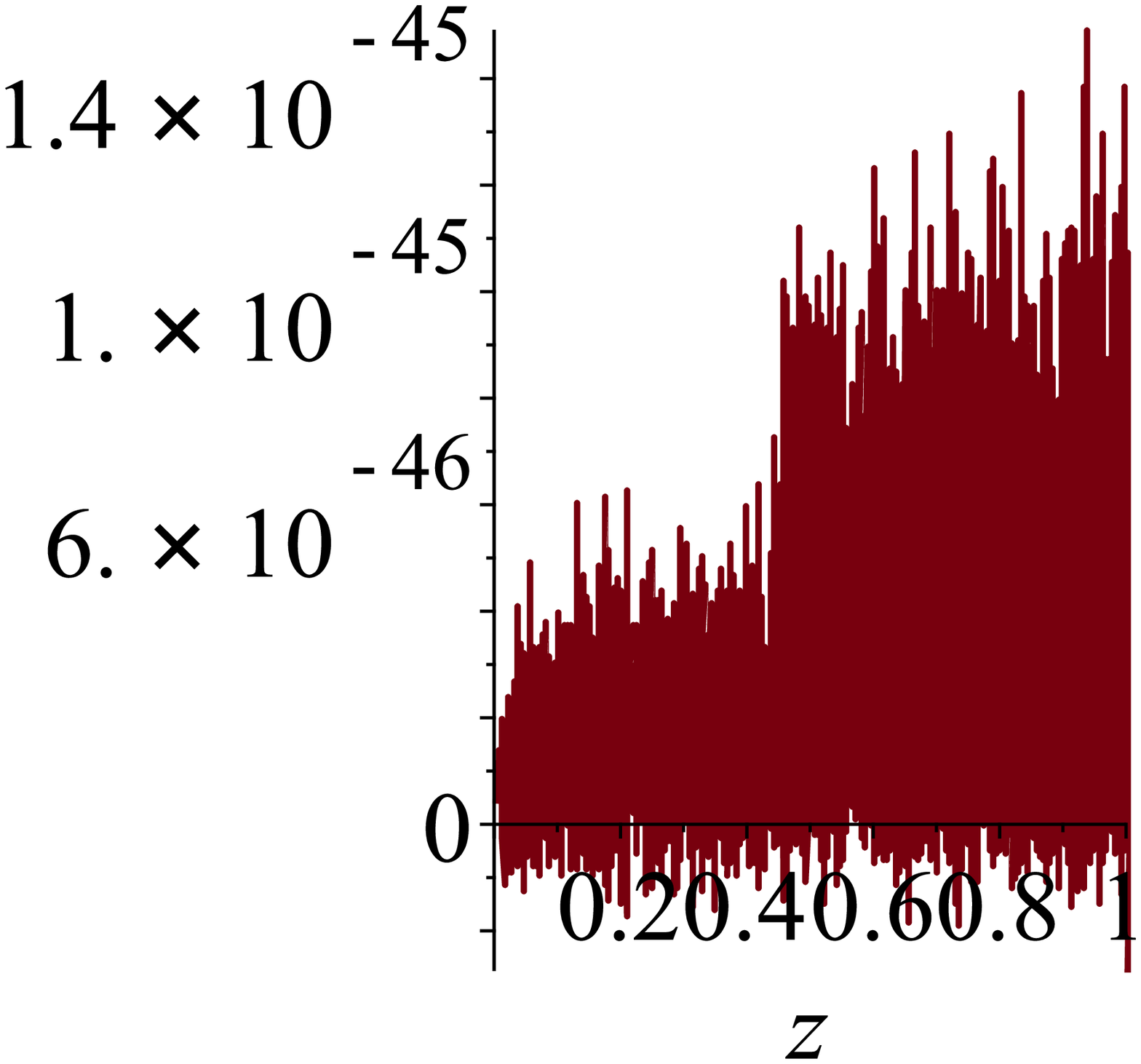}
  \includegraphics[,width=0.45\textwidth]{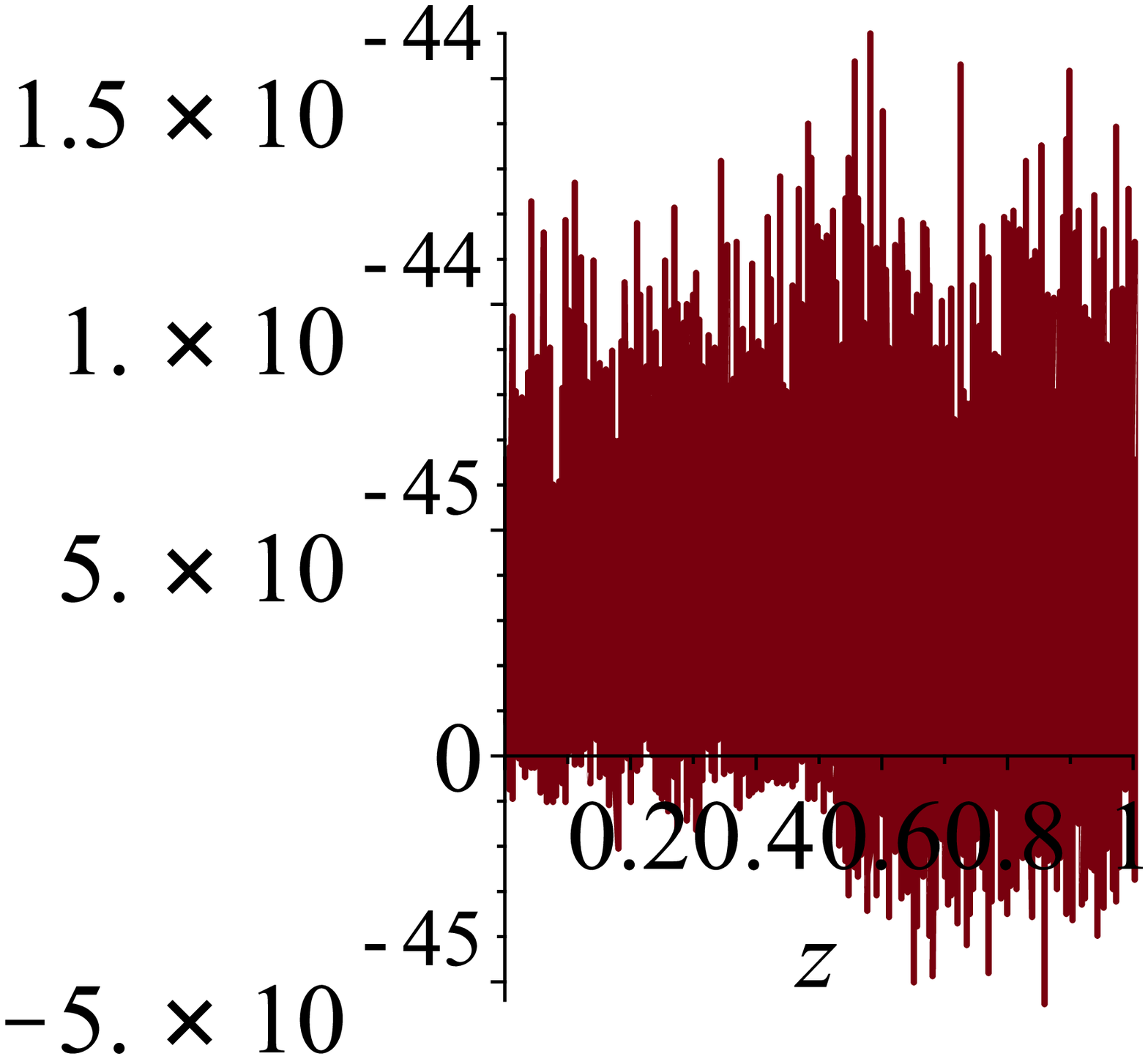}\\
 \vspace{-60pt}
Fig. 4. {\small Graphs of absolute errors for  $N=5$ and $N=10$: $|\sin(z)-P_{\sin(z)}(z)|$.}
\end{center}

The Figure 4 shows exactness our optimal interpolation formula  for the function $\sin(x)$.
\medskip

\medskip
{\small
\begin{tabular}{p{9cm}}
Bukhara state university, Bukhara, Uzbekistan\\
\end{tabular}
}

\label{Babaev2}
\end{document}